\documentstyle[11pt]{article}

\newcommand{\ff}

\newcommand{\sindex}{\end{tabular}\hspace{1em}\begin{minipage}[t]{\textwidth}
   \begin{tabbing}}

\newenvironment{proof}{\noindent{\bf Proof }}{\hfill\rule{2mm}{2mm}}
\newtheorem{theorem}{Theorem}[section]
\newtheorem{lemma}{Lemma}[section]
\newtheorem{corollary}{Corollary}[section]
\newtheorem{proposition}{Proposition}[section]

\newtheorem{defintion}{Defintion}[section]
\newtheorem{example}{Example}[section]

\newenvironment{theo}{\begin{theorem}\em}{\rm\end{theorem}}
\newenvironment{cor}{\begin{corollary}\em}{\rm\end{corollary}}
\newenvironment{lem}{\begin{lemma}\em}{\rm\end{lemma}}
\newenvironment{prop}{\begin{proposition}\em}{\rm\end{proposition}}

\newenvironment{defi}{\begin{defintion}\em}{\rm\end{defintion}}
\newenvironment{exa}{\begin{example}\em}{\rm\end{example}}

\begin{document}
\title{
   Permutations avoiding a set of patterns $T\subseteq S_3$ 
   \\ and a pattern $\tau\in S_4$} 

\author{
   Toufik Mansour \\
   Department of Mathematics \\
   Haifa University}
\maketitle

\subsection*{\centering Abstract}
{\em
In this paper we calculate the cardinality of the set $S_n(T,\tau)$ of all
permutations in $S_n$ that avoid one pattern from $S_4$ and a nonempty set 
of patterns from $S_3$. 

The main body of the paper is divided into four sections corresponding to the cases $|T|=1,2,3,$
and $|T|\geq 4$. At the end of each section we provide the tables accumulating all
the results obtained.}

\section{Introduction}
Let $[k]=\{1,\dots,k\}$ be a (totally ordered) {\em alphabet} on $k$ letters, and let
$\alpha\in [k]^m$, $\beta\in [l]^m$ with $l\leq k$. We say that $\alpha$ is  
{\em order-isomorphic} to $\beta$ if rhe following condition holds 
for all $1\leq i<j\leq n$: $\alpha_i<\alpha_j$ if and only if $\beta_i<\beta_j$.\\

We say that $\tau\in S_n$ {\em avoids} $\alpha\in S_k$ (with $k\leq m$) if there is exist 
$1\leq i_1\leq\dots\leq i_k\leq n$ such that $(\tau_{i_1},\dots,\tau_{i_k})$ is 
order-isomorphic to $\alpha=(\alpha_1,\dots,\alpha_k)$, and we say that $\tau$ 
{\em avoids} $\alpha$ if $\tau$ does not contain $\alpha$. 

The set of all permutations in $S_n$ avoiding $\alpha$ is denoted $S_n(\alpha)$. In a 
similar way, for any $A\subset \bigcup_{m=1}^{\infty} S_m$ we write $S_n(A)$ to denote 
the set of permutations in $S_n$ avoiding all the permutations in $A$.\\

The study of the sets $S_n(\alpha)$ was initiated by Knuth ~\cite{knu}, who peroved that 
$|S_n(\alpha)|=\frac{1}{n+1}{{2n}\choose n}$ for any $\alpha\in S_3$. Knuth's results 
where further extended in two directions. West ~\cite{wes} and Stankova ~\cite{sta} 
analyzed $S_n(\alpha)$ for $\alpha\in S_4$ and obtained the complete classifiction, which 
contains $3$ distinct cases. This classification, however does not give exact values 
of $S_n(\alpha)$. On the other hand, Simion and Schmidt ~\cite{sim} studied $S_n(T)$ 
for arbitrary subsets $T\subseteq S_3$ and discovered $7$ distinct cases. The study 
of $S_n(\alpha,\tau)$ for all $\alpha\in S_3$, $\tau\in S_4$, $\tau$ avoids $\alpha$, 
was completed by West ~\cite{wes}, Billey, Jockusch and Stanley ~\cite{bjs} and Guibert 
~\cite{gui}. In the persent paper we continue this work and calculate the cardinalities 
of the sets $S_n(T,\alpha)$ for all nonempty sets $T\subseteq S_3$ and all permutations 
$\alpha\in S_4$. The rest of the intorduction contains several auxilary definitions and 
results.\\

We define two symmerties on permutations, the {\em reversal} $r:S_n\rightarrow S_n$ as
follows:
  $$r:(\alpha_1,\alpha_2,\dots,\alpha_n)\longmapsto (\alpha_n,\dots,\alpha_2,\alpha_1),$$
and the {\em inverse} $i:S_n\rightarrow S_n$ as follows:

$$i:\alpha\longmapsto {\alpha}^{-1}.$$

These symmetries can be extended to an arbitrary set of permutations 
$T\subset \bigcup_{m=1}^{\infty} S_m$ as follows:
\begin{center}
	$r(T)=\{r(\alpha)\mid \alpha\in T\}$, \\
	$T^{-1}=\{\alpha^{-1}\mid \alpha\in T\}.$ 
\end{center}

We denote by $M_p$ the group of transformations of $\bigcup_{m=1}^{\infty} S_m$ generated by 
$r$ and $i$, and we define {\em Symmetry classes} as the orbits of $M_p$ in 
$\bigcup_{m=1}^{\infty} S_m$. In other words, $A,B$ belongs to be same symmetry class 
$\mathcal{T}$ if there exist $g\in M_p$ such that $g(A)=B$.

\begin{prop}
   Let $A,B\subset \bigcup_{m=1}^{\infty} S_m$ belong to be same symmetry class $\mathcal{T}$.
   Then $|S_n(A)|=|S_n(B)|$.
\label{symop}
\end{prop}
\begin{proof}
  By Burstein ~\cite{bur} and defintions.
\end{proof}
\\

Let $b_1<\dots<b_n$; we denote by $S_{\{b_1,\dots,b_n\}}$ the set of all permutation of 
the numbers $b_1,\dots,b_n$; for example, $S_{\{1,\dots,n\}}$ is just $S_n$. As above 
we denote by $S_{\{b_1,\dots,b_n\}}(T)$ the set of all permutations in 
$S_{\{b_1,\dots,b_n\}}$ avoiding all the permutations in $T$.

\begin{prop}
  Let $\tau\in S_{\{c_1,\dots,c_k\}}$. Then there exists 
  permutation $\alpha\in S_k$ such that $|S_n(\alpha)|=|S_{\{b_1,\dots,b_n\}}(\tau)|$.
\label{oneword prop}
\end{prop}
\begin{proof}
  We define a function $f:S_{\{c_1,\dots,c_k\}}\rightarrow S_k$ by
  \begin{center} $f((c_{i_1},c_{i_2},\dots,c_{i_k}))=(i_1,i_2,\dots,i_k).$ \end{center}
  Then evidently $|S_{\{b_1,\dots,b_n\}}(\tau)|=|S_n(f(\tau))|$.
\end{proof}

\begin{cor}
   Let $T\subseteq S_{\{c_1,\dots,c_k\}}$. Then 
   there exists $R\subseteq S_k$ such that $|S_n(R)|=|S_{\{b_1,\dots,b_n\}}(T)|$.
\label{setword cor}
\end{cor}
\begin{proof}
  Let $T=\{\tau_1,\dots,\tau_l\}$ and $R=\{f(\tau_1),\dots,f(\tau_l)\}$ 
  where $f$ is defined in proposition \ref{oneword prop} . On the other hand by defintions 
       $$S_{\{b_1,\dots,b_n\}}(T)=\bigcap_{i=1}^l S_{\{b_1,\dots,b_n\}}(\tau_i).$$
  Hence by the isomorphism $f$ and proposition \ref{oneword prop} this corollary holds.
\end{proof}
\section{$|S_n(\alpha,\tau)|$ for all $\alpha\in S_3$ and $\tau\in S_4$}
\begin{defi}
  Let $\tau\in S_k$. We define $A_{\tau}^m$ as the set of all 
  the words $\alpha=(\alpha_1,\dots,\alpha_k)$ such that $1\leq \alpha_i\leq m$ 
  for all $i$ and $\alpha$ is order-isomorphic to $\tau$.
\end{defi}
\begin{prop}
  $|A_{\tau}^m|={m \choose k}$ for any $\tau\in S_k$ and $m\geq k$.
\end{prop}
\begin{proof}
  If $m=k$ then the proposition is trivial. We proceed by induction. Let 
  $|A_{\tau}^m|={m \choose k}$, and we want to calculate the cardinality of the set 
  $A_{\tau}^{m+1}$. Assume that $\tau_d=k$. We have
	$$A_{\tau}^{m+1}=\bigcup_{i=k}^{m+1} A_{\tau_1,\dots,\tau_{d-1},\tau_{d+1},\dots,\tau_k}^{i-1},$$
  and the sets in above relation are disjoint. Hence by induction we obtain 
  $|A_{\tau}^{m+1}|={k-1 \choose k-1}+{k \choose k-1}+\dots+{m \choose k-1}={m+1 \choose k}$.
\end{proof}

\begin{defi} 
  Let $\tau\in S_k$. We define $V_{\tau}^m$ as the set of all the 
  permutations $\alpha\in S_m$ for which there exist $1\leq i_1<i_2<\dots<i_k\leq m$ 
  such that $(\alpha_{i_1},\dots,\alpha_{i_k})\in A_{\tau}^m$.
\end{defi}

\begin{exa}  
  $A_{132}^4=\{132,142,143,243\}$.\\ 
  $V_{132}^4=\{1324,1342,1432,4132,1423,3142,1243,2143,2431,2413\}$, and by
  direct calculations we have $|V_{\tau}^4|=10$ for all $\tau\in S_3$.
\end{exa}

\begin{defi}
  Let $\nu:P(S_k)\rightarrow P(S_{k+1})$ be defined by 
  $\nu(T)=\bigcup_{\tau\in T} V_{\tau}^{k+1}$ where $P(X)$ as usual stands for the 
  power set of $X$; we denote $\nu^p(T)=\nu(\nu^{p-1}(T))$ for $p>1$.
\end{defi}

\begin{lem}
\label{onemu lem} 
  $S_n(\nu(\tau))=S_n(\tau)$ for all $\tau\in S_k$ and $n\geq k+2$.
\end{lem}
\begin{proof}
  Let $\alpha\in S_n(\tau)$ and assume that $\alpha\notin S_n(\nu(\tau))$, 
  which means that $\alpha$ contains a subsequence $\beta\in \nu(\tau)$. So
  by defintions $\alpha$ contains a subsequence $\gamma$ such that $\gamma$ is order-isomorphic to $\tau$, 
  a contradiction. Hence $S_n(\tau)\subseteq S_n(\nu(\tau))$.

  On the other hand, let $\alpha$ contain $\tau$, which means that  there
  exist $1\leq i_1<\dots<i_k\leq n$ such that $(\alpha_{i_1},\dots,\alpha_{i_k})$ is order-isomorphic to $\tau$.
  Let $\alpha_{i_{m_1}}$ be the maximal element in $\{\alpha_{i_1},\dots,\alpha_{i_k}\}$. 
  If $\alpha_{m_1}\neq n$ then the subsequence 
  $\alpha_{i_1},\dots,\alpha_{i_{d-1}},n,\alpha_{i_{d+1}},\dots,\alpha_{i_k}$ of $\alpha$
  is order-isomorphic to some permutation in the set $V_{\tau}^{k+1}$, so 
  $S_n(\nu(\tau))\subseteq S_n(\tau)$. Hence we can assume that $\alpha_{m_1}=n$. 
  Let $\alpha_{m_2}$ be the maximal element in $\{\alpha_{i_1},\dots,\alpha_{i_k}\}\backslash\{\alpha_{m_1}\}$; 
  by the same reason we see that the only nontrivial case is $\alpha_{m_2}=n-1$ an so on. So
  $(\alpha_{i_1},\dots,\alpha_{i_k})$ is just a permutation of the numbers $n+1-k,\dots,n$.

  Now, since $n\geq k+2$, there exists $d$ such that the subsequence
  $\alpha_{i_1},\dots,\alpha_{i_{d-1}}$, $1,\alpha_{i_{d+1}},\dots,\alpha_{i_k}$
  of $\alpha$ is order-isomorphic to some permutation in the set $\nu(\tau)$. 

  Hence in any case $\alpha$ contains some permutation $\beta$ with $\beta\in \nu(\tau)$, which means that
  if $\alpha$ contains $\tau$ then $S_n(\nu(\tau))\subseteq S_n(\tau)$.
\end{proof}

\begin{cor}
\label{setmu cor}
  $S_n(\nu(T))=S_n(T)$ for all $T\subseteq S_k$ and $n\geq k+2$.
\end{cor}
\begin{proof}
  By definitions $S_n(\nu(T))=\bigcap_{\beta\in T} S_n(\nu(\beta))$, 
  and by lemma \ref{onemu lem} we obtain $S_n(\nu(T))=\bigcap_{\beta\in T} S_n(\beta)$. 
  Hence, again by definition, we obtain $S_n(\nu(T))=S_n(T)$.
\end{proof}

\begin{theo} 
\label{mu theo}
  $S_n(\nu^p(T))=S_n(T)$ for all $T\subseteq S_k$ and $n\geq k+p+1$.
\end{theo}
\begin{proof} 
  By definitions, corollary \ref{setmu cor} and induction.
\end{proof}

\begin{theo}
\label{opmu theo}
  Let $\alpha\in S_k$ and $\tau\in S_m$ with $k<m$. Then
  $\tau$ contains $\alpha$ if and only if $S_n(\alpha,\tau)=S_n(\alpha)$.
\end{theo}
\begin{proof} 
  Assume first that $\tau$ contains $\alpha$. By defintions we have $S_n(\alpha,\tau)=S_n(\alpha)\cap S_n(\tau)$. 
  On the other hand $\tau\in\nu^{m-k}(\alpha)$, and by theorem \ref{mu theo} we obtain $S_n(\alpha)=S_n(\nu^{m-k}(\alpha))$, 
  which means that $S_n(\alpha,\tau)=S_n(\nu^{m-k}(\alpha))$, 
  so again by theorem \ref{mu theo} we obtain $S_n(\alpha,\tau)=S_n(\alpha)$.

  Now let $S_n(\alpha,\tau)=S_n(\alpha)$, then $S_n(\alpha)\subseteq S_n(\tau)$, which means that 
  if $\beta\notin S_n(\tau)$ then $\beta\notin S_n(\alpha)$. By taking $\beta=\tau$, we get that 
  $\tau$ must contain $\alpha$. 
\end{proof}

\begin{cor}
\label{muknuth cor}
  Let $\alpha\in S_3$, $\tau\in S_4$ and $\tau$ contain $\alpha$. Then $|S_n(\alpha,\tau)|=c_n$ 
  where $c_n$ is the $n$-th Catalan number.
\end{cor}
\begin{proof} 
  By the theorem \ref{opmu theo} and Knuth ~\cite{knu}.
\end{proof}

{ \footnotesize
\begin{table}[h]
    \begin{tabular}{|l|l|l|l|} \hline
						&			&						&			\\
	\emph{Representative $T\in\mathcal{T}$} &\emph{$|\mathcal{T}|$} & \emph{$|S_n(T)|$ for $T\in\mathcal{T}$}  	& \emph{Reference}	\\ \hline\hline

	$\{\alpha,\tau\}$ when $\alpha\in S_3$,			&	&						&			\\
	$\tau\in S_4$ and $\tau$ contains $\alpha$		& $60$ 	& $c_n=\frac{1}{n+1} {{2n} \choose {n}}$	& corollary \ref{muknuth cor}\\ \hline

								&	&						&			\\
	$\overline{\{123,1432\}}$, $\overline{\{123,2143\}}$	&	&						&			\\
	$\overline{\{123,2413\}}$, $\overline{\{132,1234\}}$	& 	&  						&			\\ 
	$\overline{\{132,2134\}}$, $\overline{\{132,2314\}}$	& $46$	& $f_{2n-2}$ where $f_n$ is the 		& West ~\cite{wes}	\\
	$\overline{\{132,2341\}}$, $\overline{\{132,3241\}}$	&	& $n$-th Fibonacci number			&			\\
	$\overline{\{132,3412\}}$					&	&					&			\\ \hline

								&	&						&			\\
	$\overline{\{132,3421\}}$, $\overline{\{132,4231\}}$	& $12$ 	& $1+(n-1)2^{n-2}$				& West ~\cite{wes} , Guibert ~\cite{gui} resp.\\ \hline

								&	&						&			\\
	$\overline{\{123,2431\}}$				& $8$	& $3\cdot 2^{n-1}-{{n+1}\choose{2}}-1$		& West ~\cite{wes}	\\ \hline

								&	&						&			\\
	$\overline{\{123,3421\}}$				& $4$	& ${n \choose  4}+2{n\choose 3}+n$		& West ~\cite{wes}	\\ \hline

								&	&						&			\\
	$\overline{\{132,3214\}}$				& $4$	& $\frac{(1-x)^3}{1-4x+5x^2-3x^3}$		& West ~\cite{wes}	\\ \hline
	
								&	&						&			\\
	$\overline{\{132,4321\}}$				& $4$	& ${n\choose 4}+{{n+1}\choose 4}+{n\choose 2}+1$& West ~\cite{wes}	\\ \hline 

								&	&						&			\\
	$\overline{\{123,4321\}}$				& $2$	& $0$						& Erd\"os and Szekeres ~\cite{erd}\\ \hline

								&	&						&			\\
	$\overline{\{123,3412\}}$					& $2$	& $2^{n+1}-{{n+1}\choose {3}}-2n-1$	& Billey, Jockusch and Stanley ~\cite{bjs}\\ \hline

								&	&						&			\\
	$\overline{\{123,4231\}}$				& $2$	& ${n\choose 5}+2{n\choose 4}+{n\choose 3}+{n\choose 2}+1$ & West ~\cite{wes}	\\ \hline
   \end{tabular}
   \caption {Cardinalities of the sets $S_n(\alpha,\tau)$ when $\alpha\in S_3$ and $\tau\in S_4$.}
   \label{tab34}
\end{table} }
\normalsize

Billey, Jockusch and Stanley ~\cite{bjs} show that $|S_n(123,3412)|=2^{n+1}-{{n+1} \choose {3}} -2n-1$, 
Guibert ~\cite{gui} show that $|S_n(132,4231)|=1+(n-1)2^{n-2}$, Erd\"os and Szekeres ~\cite{erd} 
show that $|S_n((1,\dots,l),(k,\dots,1))|=0$. West ~\cite{wes} and corollary \ref{muknuth cor} 
compelete all the calculation of the cardinalities of $S_n(\alpha,\tau)$ when 
$\alpha\in S_3$ and $\tau\in S_4$. All these results we summarize in the above 
table. \\

The first column contains the list of representatives for the symmetry classes, 
one representative per each class. The second column contains the number of 
sets in each symmetry class, the third column contains expression for the 
cardinality of the sets in the corresponding symmetry class. The last column 
provides a reference to the paper (or theorem in the present papers) where 
this expression is obtained. 
\section{$|S_n(\alpha_1,\alpha_2,\tau)|$ when $\alpha_1\neq\alpha_2\in S_3$ and $\tau\in S_4$}
In this section we calculate the cardinalities of all the sets 
$S_n(\alpha_1,\alpha_2,\tau)$ for any two different permutations 
$\alpha_1,\alpha_2\in S_3$ and $\tau\in S_4$.
\begin{prop}
 \label{knu44}
 Let $\alpha_1$, $\alpha_2$ be two different permutations in $S_3$ and let $\tau$ 
 be permutation in $S_4$ such that $\tau$ contains $\alpha_1$ or $\alpha_2$. Then
 $|S_n(\alpha_1,\alpha_2,\tau)|=|S_n(\alpha_1,\alpha_2)|$.
\end{prop}
\begin{proof}
  By theorem \ref{opmu theo} and definitions.
\end{proof}
\\

By proposition \ref{knu44}, Simion and Schmidt ~\cite{sim} and Erd\"os and Szekeres ~\cite{erd} we obtain 
the following theorem.

\begin{theo} 
\label{the0n0}
  For all $n\in\mathcal{N}$, $\tau\in S_4$:
	\begin{enumerate}
	\item $|S_n(\alpha_1,\alpha_2,\tau)|=2^{n-1}$ if $\tau$ contains $\alpha_1$ or 
		$\alpha_2$ and $(\alpha_1,\alpha_2)=(123,132)$, $(132,213)$, 
		$(132,231)$, $(132,231)$ or $(132,312)$.

	\item $|S_n(123,321,\tau)|=0$ for all $n\geq 5$.

	\item $|S_n(\alpha_1,\alpha_2,\tau)|=0$ if $\alpha_1\in S_3$, 
		$\tau$ contains $\alpha_2$, $n\geq 7$ and $\alpha_2=123$ or 
		$\alpha_2=321$. 
	\end{enumerate} 
\end{theo}

Let us analyse other cases.

\begin{theo} 
\label{thtri}
  For all $n\in \mathcal{N}$,
     $$|S_n(123,132,3214)|=|S_n(123,213,1432)|=|S_n(132,213,1234)|=t_n,$$ 
where $t_n$ is the $n$-th Tribonacci number ~\cite{fei}.
\end{theo}
\begin{proof} 
  $1.$ Let $\alpha\in G_n=S_n(123,132,3214)$, and let us consider the possible value 
  of $\alpha_1$:
  \begin{description}
	\item[$1.1$]	$\alpha_1\leq n-2$. Evidently there exist $\alpha_{i_j}=n+1-j$, $j=1,2$, hence
		$\alpha$ either contains $(\alpha_1,\alpha_{i_1},\alpha_{i_2})$ which is order-
		isomorphic to $132$, or contains $(\alpha_1,\alpha_{i_2},\alpha_{i_1})$ which is 
		order$-$isomorphic to $123$, a contradiction.

	\item[$1.2$]	$\alpha_1=n-1$. Similarly to case $1.1$ we have $\alpha_2=n$ or $\alpha_2=n-2$. If $\alpha_2=n$ then $\alpha\in G_n$
		if and only if $(\alpha_3,\dots,\alpha_n)\in G_{n-2}$. If $\alpha_2=n-2$ then $\alpha_3=n$, since 
		otherwise $\alpha$ contains $3214$, hence $\alpha\in G_n$ if and only if 
		$(\alpha_4,\dots,\alpha_n)\in G_{n-3}$.
 
	\item[$1.3$]	$\alpha_1=n$. Evidently $\alpha\in G_n$ if and only if 
		$(\alpha_2,\dots,\alpha_n)\in G_{n-1}$.
  \end{description}
  Since the above cases are disjoint we obtain $|G_n|=|G_{n-1}|+|G_{n-2}|+|G_{n-3}|$.\\

$2.$ Let $\alpha\in G_n=S_n(123,213,1432)$, and let us consider the possible 
value of $\alpha_1$:
\begin{description}
	\item[$2.1$]	$\alpha_1\leq n-3$. Evidently there exist $\alpha_{i_j}=n+1-j$, 
		$j=1,2,3$. Since $\alpha$ avoids $123$ we get 
		$\alpha_{i_1}<\alpha_{i_2}<\alpha_{i_3}$, hence $\alpha$ contains 
		$(\alpha_1,\alpha_{i_1},\alpha_{i_2},\alpha_{i_3})$ which is order$-$
		isomorphic to $1432$, a contradiction.

	\item[$2.2$]	$\alpha_1=n-2$. Let $\alpha_{i_j}=n+1-j$, $j=1,2$. Since 
		$\alpha$ avoids $123$ we get $i_1<i_2$. If $i_1\geq 3$ or $i_2\geq 4$ 
		then $\alpha$ contains $213$, a contradiction. So 
		$\alpha=(n-2,n,n-1,\alpha_4,\dots,\alpha_n)$. Hence $\alpha\in G_n$ if 
		and only if $(\alpha_4,\dots,\alpha_n)\in G_{n-3}$.

	\item[$2.3$]	$\alpha_1=n-1$. If $\alpha_2\leq n-2$ then $\alpha$ contains $213$, so we have 
		$\alpha=(n-1,n,\alpha_3,\dots,\alpha_n)$. Hence $\alpha\in G_n$
		if and only if $(\alpha_3,\dots,\alpha_n)\in G_{n-2}$.
 
	\item[$2.4$]	$\alpha_1=n$. Evidently $\alpha\in G_n$ if and only if 
		$(\alpha_2,\dots,\alpha_n)\in G_{n-1}$.
\end{description}
Since the above cases are disjoint we obtain $|G_n|=|G_{n-1}|+|G_{n-2}|+|G_{n-3}|$.\\

$3.$ Let $\alpha\in G_n=S_n(132,213,1234)$, and let us consider the possible 
value of $\alpha_1$:
\begin{description}
	\item[$3.1$]	$\alpha_1\leq n-3$. Since $\alpha$ avoids $132$ we get, similarly to case $2.1$,  
		that $\alpha$ contains $(\alpha_1,n-2,n-1,n)$, 
		which is order$-$isomorphic to $1234$, a contradiction. 

	\item[$3.2$]	$\alpha_1=n-2$. If $\alpha_2\leq n-2$ then $\alpha$ contains 
		$213$, and if $\alpha_2=n$ then $\alpha$ contains $132$, so we have 
		that $\alpha_2=n-1$. If $\alpha_3\leq n-2$ then $\alpha$ contains 
		$213$, so $\alpha_3=n$. Hence $\alpha=(n-2,n-1,n,\alpha_4,\dots,\alpha_n)$. 
		That is, $\alpha\in G_n$ if and only if $(\alpha_4,\dots,\alpha_n)\in G_{n-3}$.

	\item[$3.3$]	$\alpha_1=n-1$. Similarly to case $2.3$ we have $\alpha_2=n$, 
		so $\alpha=(n-1,n,\alpha_3,\dots,\alpha_n)$. 
		Hence $\alpha\in G_n$ if and only if $(\alpha_3,\dots,\alpha_n)\in G_{n-2}$.
 
	\item[$3.4$]	$\alpha_1=n$. Evidently $\alpha\in G_n$ if and only if 
		$(\alpha_2,\dots,\alpha_n)\in G_{n-1}$.
\end{description}
Since the above cases are disjoint we obtain $|G_n|=|G_{n-1}|+|G_{n-2}|+|G_{n-3}|$.
\end{proof}

\begin{theo}  
\label{thefn2}                 
  For all $n\in\mathcal{N}$,
		$$|S_n(123,132,3241)|=|S_n(132,213,2341)|=f_{n+2}-1,$$
  where $f_n$ is the $n$-th Fibonacci number. 
\end{theo}
\begin{proof} $1.$ Let $\alpha\in G_n=S_n(123,132,3241)$ and let us consider the possible 
 value of $\alpha_1$:
 \begin{description}
	\item[$1.1$]	$\alpha_1\leq n-2$. Impossible, similarly to case $1.1$ in theorem \ref{thtri}.

	\item[$1.2$]	$\alpha_1=n-1$. Similarly to case $1.2$ in theorem \ref{thtri} we have 
		that $\alpha_2=n-2$ or $\alpha_2=n$. If $\alpha_2=n-2$ then
		we get $\alpha_n=n$ since otherwise $\alpha$ contains $3241$. Besides, 
		$\alpha_i>\alpha_j$ for all $3\leq i<j<n$ since otherwise $\alpha$ contains 
		$123$. So there is only one such permutation, namely $(n-1,n-2,\dots,1,n)$. 
		If $\alpha_2=n$ then $\alpha\in G_n$ if and only if $(\alpha_3,\dots,\alpha_n)\in G_{n-2}$.

	\item[$1.3$]	$\alpha_1=n$. Evidently $\alpha\in G_n$ if and only if $(\alpha_2,\dots,\alpha_n)\in G_{n-1}$.
\end{description}
Since the above cases are disjoint we obtain $|G_n|=|G_{n-1}|+|G_{n-2}|+1$.
By the transformation $g_n=|G_n|+1$ we get $g_n=g_{n-1}+g_{n-2}$. Besides $g_4=8$ and $g_5=13$, which means 
that $g_n=f_{n+2}$, hence $|G_n|=f_{n+2}-1$ for all $n\geq 4$. It is easy to see that for $n=1,2,3$ the same
formula holds.\\

$2.$ Let $\alpha\in K_n=S_n(132,213,2341)$, and let us consider the possible value of 
$\alpha_1$:
\begin{description}
	\item[$2.1$] 	$2\leq\alpha_1\leq n-2$. If $\alpha_2=n$ then $\alpha$ contains $132$. Let
		$\alpha_2=n-1$, if $\alpha_3=n$ then $\alpha$ contains $2341$ otherwise $\alpha$ contains $213$. 
		Hence $\alpha_2\leq n-2$. If $\alpha_2<\alpha_1$ then $\alpha$ contains $213$, so we have that $\alpha_2>\alpha_1$. 
		If $\alpha_3>\alpha_2$ then $\alpha$ contains $2341$, otherwise $\alpha$ contains $213$, a contradiction.

	\item[$2.2$] 	$\alpha_1=1$. Since $\alpha$ avoids $132$, we have only one permutation $(1,\dots,n)$.

	\item[$2.3$]	$\alpha_1=n-1$. Since $\alpha$ avoids $213$, we have that $\alpha_2=n$. Hence 
	      $\alpha\in K_n$ if and only if $(\alpha_3,\dots,\alpha_n)\in K_{n-2}$.

	\item[$2.4$] 	$\alpha_1=n$. Evidently $\alpha\in K_n$ if and only if $(\alpha_2,\dots,\alpha_n)\in K_{n-1}$.
\end{description}
Since the above cases are disjoint we obtain $|K_n|=|K_{n-1}|+|K_{n-2}|+1$; besides 
$|K_4|=7$ and $|K_5|=12$. Similarly to the first part of the proof we obtain 
$|K_n|=|G_n|=f_{n+2}-1$ for all $n\in\mathcal{N}$.
\end{proof}

\begin{theo} 
\label{the3n5}
  For all $n\geq 3$,
		$$|S_n(123,132,3421)|=|S_n(123,213,3421)|=3n-5.$$ 
\end{theo}
\begin{proof} $1.$ Let $\alpha\in G_n=S_n(123,132,3421)$, and let us consider the 
possible value of $\alpha_1$:
\begin{description}
	\item[$1.1$] 	$\alpha_1\leq n-2$. Impossible, similarly to case $1.1$ in theorem \ref{thtri}.

	\item[$1.2$] 	$\alpha_1=n-1$. Consider $k$ such that $\alpha_k=n$. Since $\alpha$ 
		avoids $3421$ we have $\alpha_i<\alpha_j$ for all $k<i<j\leq n$, 
		hence if $k\leq n-3$ then $\alpha$ contains $123$. So either $k=n$ or 
		$k=n-1$ or $k=n-2$. Besides $\alpha_i>\alpha_j$ for all $1<i<j<k$ since 
		otherwise $\alpha$ contains $123$, and $\alpha_j>\alpha_i$ for 
		$i<k<j$ since otherwise $\alpha$ contains $132$. Hence in this case there are three possible 
		permutations: $(n-1,\dots,1,n)$, $(n-1,\dots,2,n,1)$ and $(n-1,\dots,3,n,1,2)$.

	\item[$1.3$] 	$\alpha_1=n$. Evidently $\alpha\in G_n$ if and only if 
		$(\alpha_2,\dots,\alpha_n)\in G_{n-1}$.
\end{description}
Since the above cases are disjoint we obtain $|G_n|=|G_{n-1}|+3$. Besides $|G_3|=4$, hence
$|G_n|=3n-5$. \\

$2.$ Let $\alpha\in K_n=S_n(123,213,3421)$, and let us consider the possible value of 
$\alpha_1$:
\begin{description}
	\item[$2.1$] 	$4\leq\alpha_1\leq n-1$. Let $\alpha_{i_j}=j$, $j=1,2,3$, 
		and $\alpha_{i_4}=n$. Since $\alpha$ avoids $213$ we have that $i_4<i_j$ 
		for all $j=1,2,3$. On the other hand $\alpha$ avoids $3421$, so we have 
		$i_1<i_2<i_3$, which means that $\alpha$ contains $123$, a contradiction.

	\item[$2.2$] 	$\alpha_1=n$. Evidently $\alpha\in K_n$ if and only if 
		$(\alpha_2,\dots,\alpha_n)\in K_{n-1}$.

	\item[$2.3$] 	$\alpha_1=1$. Since $\alpha$ avoids $123$ we have only one 
		permutation $(1,n,\dots,2)$.

	\item[$2.4$] 	$\alpha_1=2$. Since $\alpha$ avoids $123$ we have that $\alpha$ contains $(n,\dots,3)$. Since 
		$\alpha$ avoids $213$ we have $\alpha_n=1$. So we have only one permutation $(2,n,\dots,3,1)$.

	\item[$2.5$] 	$\alpha_1=3$. Since $\alpha$ avoids $123$ and $213$ we have two permutations 
		$(3,n,\dots,4,1,2)$ or $(3,n,\dots,4,2,1)$, but $\alpha$ avoids $3421$, so we
		we have only one permutation $(3,n,\dots,4,1,2)$.
\end{description}
Since the above cases are disjoint we obtain $|K_n|=|K_{n-1}|+3$. Besides $|K_3|=4$, hence
$|K_n|=|G_n|=3n-5$.
\end{proof}

\begin{theo} 
\label{the2n2}
Let $\alpha\in \{1432,2143,2431,3214,3241,3421\}$, 
then for all $n\geq 2$, 
		$$|S_n(123,312,\alpha)|=2n-2.$$
\end{theo}
\begin{proof} $1.$ Let $\alpha\in G_n=S_n(123,312,1432)$. Fix $k$ such that $\alpha_k=1$ 
and let us consider the possible value of $k$:
\begin{description}
	\item[$1.1$]	$k=n$. Evidently $\alpha\in G_n$ if and only if 
		$(\alpha_1,\dots,\alpha_{n-1})\in S_{\{2,\dots,n\}}(123,312,1432)$.

	\item[$1.2$]	$k=n-1$. If $\alpha_n<n$ then $\alpha$ contains $312$, so $\alpha_n=n$. On the other hand
		$\alpha$ avoids $123$, so we have only one permutation $(n-1,\dots,1,n)$.

	\item[$1.3$]	$k=n-2$. Similarly to the above case, we have $\alpha=(n-2,\dots,1,n,n-1)$ or
		$\alpha=(n-2,\dots,1,n-1,n)$, but the permutation $(n-2,\dots,1,n-1,n)$ contains $123$, 
		so we have only one permutation $(n-2,\dots,n,n-1)$.

	\item[$1.4$]	$1\leq k\leq n-3$. Since $\alpha$ avoids $123$ we have that $\alpha_{n-2}>\alpha_{n-1}>\alpha_n$ and
		$\alpha$ contains $1432$, a contradiction.
		
\end{description}
Since the above cases are disjoint we obtain by corollary \ref{setword cor} that 
$|G_n|=|G_{n-1}|+2$. Besides $|G_3|=4$, hence $|G_n|=2n-2$. \\

$2.$ Let $\alpha\in G_n=S_n(123,312,2143)$. Fix $k$ such that $\alpha_k=1$ 
and let us consider the possible value of $k$:
\begin{description}
	\item[$2.1$]	$k=n$ or $k=n-1$. Similarly to cases $1.1$ or $1.2$ resp.

	\item[$2.2$]	$k=1$. Since $\alpha$ avoids $123$ we have only one permutation $(1,n,\dots,2)$.

	\item[$2.3$]	$2\leq k\leq n-2$. Since $\alpha$ avoids $123$ we have that $\alpha_{n-1}>\alpha_n$; 
		since $\alpha$ avoids $312$ we have $\alpha_n>\alpha_i$ and $\alpha_{n-1}>\alpha_i$ for all $i<k$, so
		$\alpha$ contains $2143$, a contradiction.		
\end{description}
Since the above cases are disjoint we obtain $|G_n|=|G_{n-1}|+2$. Besides $|G_3|=4$, 
hence $|G_n|=2n-2$.\\

$3.$ Let $A=\{123,312,2431\}$ and $G_n=S_n(r(A^{-1}))=S_n(321,132,2314)$. By 
proposition $1.1$ we get $|G_n|=|S_n(A)|$. Let $\alpha\in G_n$ and $\alpha_1=t$. 
If $\alpha_2\geq t+2$ then $\alpha$ contains $132$, and if $2\leq\alpha_2\leq t-1$ 
then $\alpha$ contains $321$, so we have $\alpha_2=1$ or $\alpha_2=t+1$.\\

Let $\alpha_2=1$; since $\alpha$ avoids $132$ we have that 
$\alpha=(t,1,\dots,t-1,t+1,\dots,n)$, anthere are $n-1$ permutations of this type.\\

Let $\alpha_2=t+1$; since $\alpha$ avoids $132$ we have that $\alpha$ conatins 
$(t,t+1,\dots,n)$. If $\alpha_i<t$ and $3\leq i\leq n-t+1$ then $\alpha$ contains $2314$, 
which means that $\alpha=(t,\dots,n,\alpha_{n-t+2},\dots,\alpha_n)$. Since $\alpha$ 
avoids $321$ we have that $\alpha=(t,\dots,n,1,\dots,t-1)$, 
and there are $n-1$ permutations of this type. Hence $|G_n|=2n-2$.\\
 
$4.$ Let $\alpha\in G_n=S_n(123,312,3214)$ and $\alpha_1=t$. If 
$t+1\leq\alpha_2\leq n-1$ then $\alpha$ contains $123$, and if $\alpha_2\leq t-2$ 
then $\alpha$ contains $312$, so we have that $\alpha_2=t-1$ or $\alpha_2=n$. \\

Let $\alpha_2=n$; since $\alpha$ avoids $312$ we get $\alpha=(t,n,\dots,t+1,t-1,\dots,1)$, 
and there are $n-1$ permutations of this type.\\

Let $\alpha_2=t-1$; if $t+1\leq\alpha_3\leq n-1$ then $\alpha$ contains $123$, and if 
$\alpha_3\leq t-2$ then $\alpha$ contains $3214$, so we have that $\alpha_3=n$. Since 
$\alpha$ avoids $312$ we have that $\alpha=(t,t-1,n,\dots,t+1,t-2,\dots,1)$, and there
are $n-1$ permutations of this type. Hence $|G_n|=2n-2$.\\

$5.$ Let $\alpha\in G_n=S_n(123,312,3241)$ and $\alpha_1=t$. Similarly 
to case $4$ we have that $\alpha_2=t-1$ or $\alpha_2=n$.\\

Let $\alpha_2=n$; since $\alpha$ avoids $312$ we have that $\alpha=(t,n,\dots,t+1,t-1,\dots,1)$, 
and there are $n-1$ permutations of this type.\\

Let $\alpha_2=t-1$. If $t=2$ then, since $\alpha$ avoids $123$ we get 
$\alpha=(2,1,n,\dots,3)$. If $t\geq 3$ then, since $\alpha$ avoids $312$ we have that 
$\alpha$ contains $(t-1,t-2,\dots,1)$. If there exist $3\leq i<t$ such that $\alpha_i>t$, 
then $\alpha$ contains $3241$, so $\alpha=(t,\dots,1,\alpha_{t+1},\dots,\alpha_n)$. 
Since $\alpha$ avoids $123$ we have that $\alpha=(t,\dots,1,n,\dots,t+1)$, and there are 
$n-2$ permutations of this type. Hence $|G_n|=2n-2$.\\

$6.$ Let $\alpha\in G_n=S_n(123,312,3421)$, fix $k$ such that $\alpha_k=1$ 
and let us consider the possible value of $k$:
\begin{itemize}	
	\item[$6.1$]	$k\leq n-1$. Since $\alpha$ avoids $312$ we have 
		that $\alpha_i<\alpha_j$ for all $i<k<j$, so since $\alpha$ avoids $123$ 
		we get $\alpha_i>\alpha_j$ for all $i<j<k$ or $k<i<j$. Hence
		$\alpha=(k,k-1,\dots,1,n,n-1,\dots,k+1)$, and there are $n-1$ 
		permutations of this type.

	\item[$6.2$]	$k=n$. Evidently $\alpha\in G_n$ if and only if $(\alpha_1,\dots,\alpha_{n-1})\in S_{n-1}(123,312,231)$. 
		By Simion and Schmidt ~\cite{sim} we have $|S_{n-1}(123,312,231)|=n-1$.
\end{itemize}
Since the above cases are disjoint we obtain $|G_n|=n-1+n-1=2n-2$.
\end{proof}

\begin{theo} 
\label{thenn2}
  For all $n\in\mathcal{N}$,
	         $$|S_n(\alpha_1,\alpha_2,\tau)|={n\choose 2}+1,$$
in the following cases:
	\begin{enumerate}
	\item	$\alpha_1=123$, $\alpha_2=231$ and $\tau\in S_4$ contains $123$ or $231$.

	\item	$\alpha_1=123$, $\alpha_2\in\{132,213\}$ and $\tau\in\{3412,4231\}$.	
	
	\item 	$(\alpha_1,\alpha_2,\tau)=(132,213,3412)$.

	\item 	$\alpha_1=132$, $\alpha_2=231$ and $\tau\in\{1234,2134,3124,3214\}$.

	\item 	$(\alpha_1,\alpha_2,\tau)=(213,312,3412)$.

	\item 	$\alpha_1=213$, $\alpha_2=321$ and $\tau\in\{1324,2314,1324\}$.

	\item	$(\alpha_1,\alpha_2,\tau)=(213,132,4321)$.
	\end{enumerate}
\end{theo}
\begin{proof}
$1.$ By Simion and Schmidt ~\cite{sim} we get $|S_n(123,231)|={n\choose 2}+1$, 
hence by theorem \ref{opmu theo} we have that $|S_n(123,231,\tau)|={n\choose 2}+1$ 
for all $\tau\in S_4$ containing $123$ or $231$.\\

$2.$ Let $\alpha\in G_n=S_n(123,132,3412)$ and let us consider the possible 
value of $\alpha_1$:
\begin{description}
	\item[$2.1$]	$\alpha_1\leq n-2$. Impossible, similarly to case $1.1$ in 
		theorem \ref{the3n5}

	\item[$2.2$]	$\alpha_1=n-1$. Let $\alpha_k=n$; since $\alpha$ avoids $123$ we 
		have for all $i_1<i_2<k$, $\alpha_{i_1}>\alpha_{i_2}$. Since $\alpha$ 
		avoids $3412$ we have for all $k<i_1<i_2$, $\alpha_{i_1}>\alpha_{i_2}$.
		Since $\alpha$ avoids 132 we have $\alpha=(n-1,\dots,n-k+1,n,n-k,\dots,1)$, 
		and there are $n-1$ permutations of this type.
		
	\item[$2.3$]	$\alpha_1=n$. Evidently $\alpha\in G_n$ if and only if 
		$(\alpha_2,\dots,\alpha_n)\in G_{n-1}$.
\end{description}
Since the above casses are disjoint we obtain $|G_n|=|G_{n-1}|+n-1$. Besides $|G_4|=7$, 
hence $|G_n|={n\choose 2}+1$ for all $n\geq 4$. It is easy to see that for $n=1,2,3$ the 
same formula holds for all $n\in\mathcal{N}$.\\

$3.$ Let $\alpha\in G_n=S_n(123,132,4231)$ and let us consider the possible 
value of $\alpha_1$:
\begin{description}
	\item[$3.1$]	$\alpha_1\leq n-2$. Impossible, similarly to case $1.1$ in 
		theorem \ref{the3n5}.

	\item[$3.2$]	$\alpha_1=n-1$. Consider $k$ such that $\alpha_k=n$ and $k\leq n-1$. Since
		$\alpha$ avoids $123$ we get $\alpha_i>\alpha_j$ for all $i<j<k$, and
		since $\alpha$ avoids $132$ we get $\alpha_i>\alpha_j$ for all $i<k<j$, 
		So $\alpha=(n-1,n-2,\dots,n-k+1,n,\alpha_{k+1},\dots,\alpha_n)$. Hence
		$\alpha\in G_n$ if and only if $(\alpha_{k+1},\dots,\alpha_n)\in S_{n-k}(123,132,231)$.
		By Simion and Schmidt ~\cite{sim} we get $|S_{n-k}(123,132,231)|=n-k$. \\

		If $k=n$ then since $\alpha$ avoids $123$ we have only one permutation 
		$(n-1,\dots,1,n)$.

	\item[$3.3$]	$\alpha_1=n$. Evidently $\alpha\in G_n$ if and only if 
		$(\alpha_2,\dots,\alpha_n)\in S_{n-1}(123,132,231)$.
		By Simion and Schmidt ~\cite{sim} we have $|S_{n-1}(123,132,231)|=n-1$.
\end{description}
Since the above cases are disjoint we obtain $|G_n|=1+(1+\dots+n-2)+n-1={n\choose 2}+1$.\\

$4.$ Let $\alpha\in G_n=S_n(123,213,3412)$. Let us consider the possible 
value of $\alpha_1$:
\begin{description}
	\item[$4.1$]	$\alpha_1=1$. Since $\alpha$ avoids $123$ we have that $\alpha=(1,n,\dots,2)$.

	\item[$4.2$]	$\alpha_1=t$, $2\leq t\leq n-1$. If $\alpha_2\leq t-1$ then $\alpha$ contains 
		$213$, and if $t+1\leq \alpha_2\leq n-1$ then $\alpha$ contains $123$, So $\alpha_2=n$.
		Since $\alpha$ avoids $3412$ we have that $\alpha$ contains $(t-1,\dots,1)$. \\

		Fix $k$ such that $\alpha_k=t-1$. If there $i>k$ such that $\alpha_i>t$ then
		$\alpha$ contains $(t,t-1,\alpha_i)$ which is order-isomorphic to $213$, so 
		$\alpha=(t,n,\alpha_3,\dots,\alpha_{n-t+1},t-1,\dots,1)$. Since $\alpha$ avoids 
		$123$ we have that $\alpha=(t,n,\dots,t+1,t-1,\dots,1)$, and there are $n-2$
		permutations of this type.
		
	\item[$4.2$]	$\alpha_1=n$. Evidently $\alpha\in G_n$ if and only if $(\alpha_2,\dots,\alpha_n)\in G_{n-1}$.
\end{description}
Since the above cases are disjoint we obtain $|G_n|=|G_{n-1}|+n-1$. Besides $|G_4|=7$, 
hence similarly to the second proof we get $|G_n|={n\choose 2}+1$.\\

$5.$ Let $A=\{123,213,4231\}$ and $G_n=S_n(r(A))=S_n(321,312,1324)$. By 
proposition \ref{symop} we get $|G_n|=|S_n(A)|$. Let $\alpha\in G_n$ and let us consider 
the possible value of $\alpha_1$:
\begin{description}
	\item[$5.1$]	$3\leq \alpha_1$. Since $\alpha$ avoids $312$ we have that 
		$\alpha$ contains $321$, a contradiction.

	\item[$5.2$]	$\alpha_1=2$. Consider $k$ such that $\alpha_k=1$. Since $\alpha$ 
		avoids $321$ we have $\alpha_i<\alpha_j$ for all $2\leq i<j\leq k-1$, and
		since $\alpha$ avoids $312$ we have that $\alpha_i<\alpha_j$ for all 
		$i<k<j$. Therefore $\alpha=(2,\dots,k,1,\alpha_{k+1},\dots,\alpha_n)$. Hence
		$\alpha\in G_n$ if and only if 
		$(\alpha_{k+1},\dots,\alpha_n)\in S_{\{k+1,\dots,n\}}(321,312,213)$. By
		Simion and Schmidt ~\cite{sim} and corollary $1.3$ we get that  
		$S_{\{k+1,\dots,n\}}(321,312,213)$ contains $n-k$ permutations for all 
		$k\leq n-1$, and exactly one permutation for $k=n$.

	\item[$5.3$]	$\alpha_1=1$. Similarly to case 5.2, $\alpha\in G_n$ if and only 
			if $(\alpha_2,\dots,\alpha_n)\in S_{\{2,\dots,n\}}(321,312,213)$.
			By Simion and Schmidt ~\cite{sim} and corollary \ref{setword cor} we have that $n-1$ 
			permutations.
\end{description}
Since the above cases are disjoint we obtain $|G_n|=(n-2+\dots+1)+1+n-1={n\choose 2}+1$.\\

$6.$ Let $\alpha\in G_n=S_n(132,213,3412)$ and let us consider the possible 
value of $\alpha_1$:
\begin{description}
	\item[$6.1$]	$\alpha_1=t\leq n-1$. Fix $k$ such that $\alpha_k=n$. If 
		$\alpha_i<t$ and $i<k$ then $\alpha$ contains $213$, and if $\alpha_i>t$ 
		and $i>k$ then $\alpha$ contains $132$. Since $\alpha$ avoids $132$ we 
		have that $\alpha=(t,t+1,\dots,n,\alpha_{n-t+2},\dots.\alpha_n)$, 
		and since $\alpha$ avoids $3412$ we get 
		$\alpha=(t,t+1,\dots,n,t-1,\dots,1)$, and there are $n-1$ permutations 
		of this type.

	\item[$6.2$]	$\alpha_1=n$. Evidently $\alpha\in G_n$ if and only if 
		$(\alpha_2,\dots,\alpha_n)\in G_{n-1}$.
\end{description}
Since the above cases are disjoint we obtain $|G_n|=|G_{n-1}|+n-1$. Besides $|G_4|=7$, 
hence similarly to the second proof we get $|G_n|={n\choose 2}+1$.\\

$7.$ Let $\alpha\in G_n=S_n(132,231,1234)$ and let us consider the possible value of $\alpha_1$:
\begin{description}
	\item[$7.1$] 	$\alpha_1=1$. Since $\alpha$ avoids $132$ we get $\alpha=(1,\dots,n)$, so $\alpha$ contains $1234$, 
		a contradiction.

	\item[$7.2$]	$2\leq\alpha_1\leq n-3$. Since $\alpha$ avoids $132$ the permutation $\alpha$ contains $(\alpha_1,n-2,n-1,n)$
		which is order-isomorphic to $1234$, a contradiction.

	\item[$7.3$]	$\alpha_1=n-2$. Since $\alpha$ avoids $231$ we have that 
		$\alpha_n,\alpha_{n-1}\in\{n-1,n\}$, and since $\alpha$ avoids $132$
		we get $\alpha_{n-1}=n-1$ and $\alpha_{n}=n$. Since $\alpha$ avoids
		$1234$ we have that $\alpha_i>\alpha_j$ for all $2\leq i<j\leq n-2$.
		Hence $\alpha=(n-2,n-3,\dots,1,n-1,n)$.

	\item[$7.4$]	$\alpha_1=n-1$. If $\alpha_n\leq n-2$ then $\alpha$ contains $231$,
		so $\alpha_n=n$. Evidently $\alpha\in G_n$ if and only if 
		$(\alpha_2,\dots,\alpha_{n-1})\in S_{n-2}(132,231,123)$. By Simion and 
		Schmidt ~\cite{sim} we have $|S_{n-2}(132,231,123)|=n-2$.

	\item[$7.5$]	$\alpha_1=n$. Evidently $\alpha\in G_n$ if and only if 
		$(\alpha_2,\dots,\alpha_n)\in G_{n-1}$.
\end{description}
Since the above cases are disjoint we obtain $|G_n|=|G_{n-1}|+n-2+1$. Besides $|G_4|=7$, 
hence similarly to the second proof we get $|G_n|={n\choose 2}+1$.\\

$8.$ Let $\alpha\in G_n=S_n(132,231,2134)$ and let us consider the 
possible value of $\alpha_1$:
\begin{description}
	\item[$8.1$] 	$\alpha_1=1$. Since $\alpha$ avoids $132$ we get 
		$\alpha=(1,2,\dots,n)$.

	\item[$8.2$]	$2\leq\alpha_1\leq n-2$. Since $\alpha$ avoids $132$ we get that 
		$\alpha$ contains $(\alpha_1,n-1,n)$ and since $\alpha$ avoids $231$ we 
		get that $\alpha$ contains $(\alpha_1,1,n-1,n)$ which is order-isomorphic 
		to $2134$, a contradiction.

	\item[$8.3$]	$\alpha_1=n-1$. Similarly to case $7.4$ we get $\alpha_n=n$. So 
		$\alpha\in G_n$ if and only if 
		$(\alpha_2,\dots,\alpha_{n-1})\in S_{n-2}(132,231,213)$. By Simion and 
		Schmidt ~\cite{sim} we have $|S_{n-2}(132,231,213)|=n-2$.

	\item[$8.4$]	$\alpha_1=n$. Evidently $\alpha\in G_n$ if and only if 
		$(\alpha_2,\dots,\alpha_n)\in G_{n-1}$.
\end{description}
Since the above cases are disjoint we obtain $|G_n|=|G_{n-1}|+n-2+1$. Besides $|G_4|=7$, 
hence similarly to the second proof we get $|G_n|={n\choose 2}+1$.\\

$9.$ Let $\alpha\in G_n=S_n(132,231,3124)$ and let us consider the 
possible value of $\alpha_1$:
\begin{description}
	\item[$9.1$] 	$\alpha_1=1$. Similarly to case $8.1$ we have only one 
		permutation $(1,2,\dots,n)$.

	\item[$9.2$]	$\alpha_1=2$. Since $\alpha$ avoids $231$ we have that 
		$\alpha_2=1$, and since $\alpha$ avoids $132$ we get 
		$\alpha=(2,1,3,\dots,n)$.
		
	\item[$9.3$]	$\alpha_1=t$, $3\leq t\leq n-1$. Fix $k$ such that $\alpha_k=1$. 
		Since $\alpha$ avoids $231$ we have that $\alpha_i>\alpha_j$ for all 
		$i<j\leq k$ and $\alpha_i\leq t$ for all $i\leq k$. Since $\alpha$ avoids 
		$132$ we have that $\alpha_i<\alpha_j$ for all $k<i<j$. If $\alpha_i<t$
		and $i>k$ then $\alpha$ contains either $(t,1,\alpha_i,n)$ or $(t,1,n,\alpha_i)$, 
		which means that $\alpha$ contains $3124$ or $132$, a contradiction.
		Hence $\alpha=(t,t-1,\dots,1,t+1,t+2,\dots,n)$, and there are $n-3$ permutations 
		of this type.

	\item[$9.4$]	$\alpha_1=n$. Evidently $\alpha\in G_n$ if and only if 
		$(\alpha_2,\dots,\alpha_n)\in G_{n-1}$.
\end{description}
Since the above cases are disjoint we obtain $|G_n|=|G_{n-1}|+n-3+1+1$. Besides $|G_4|=7$,
hence similarly to the second proof we get $|G_n|={n\choose 2}+1$.\\

$10.$ Let $\alpha\in G_n=S_n(132,231,3214)$ and let us consider the 
possible value of $\alpha_1$:
\begin{description}
	\item[$10.1$] 	$\alpha_1=1$. Similarly to case $8.1$ we have only one permutation 
		$(1,2,\dots,n)$.

	\item[$10.2$]	$\alpha_1=t$, $2\leq t\leq n-1$. Since $\alpha$ avoids $132$ 
		we have that $\alpha$ contains $(t,t+1,\dots,n)$. Fix $k$ such that 
		$\alpha_k=t+1$, if $\alpha_i<t$ and $i>k$ then $\alpha$ contains $231$.
		Hence $\alpha=(t,\alpha_2,\dots,\alpha_{t},t+1,\dots,n)$. Since $\alpha$ 
		avoids $3214$ we have that $\alpha=(t,1,\dots,t-1,t+1,\dots,n)$, and there 
		are $n-2$ permutations of this type.

	\item[$10.3$]	$\alpha_1=n$. Evidently $\alpha\in G_n$ if and only if 
		$(\alpha_2,\dots,\alpha_n)\in G_{n-1}$.
\end{description}
Since th above cases are disjoint we obtain $|G_n|=|G_{n-1}|+n-1$. Besides $|G_4|=7$,
hence similarly to the second proof we get $|G_n|={n\choose 2}+1$.\\

$11.$ Let $\alpha\in G_n=S_n(213,312,2341)$ and let us consider the 
possible value of $\alpha_1$:
\begin{description}
	\item[$11.1$]	$\alpha_1=1$. Evidently $\alpha\in G_n$ if and only if 
		$(\alpha_2,\dots,\alpha_n)\in S_{\{2,\dots,n\}}(213,312,2341)$.

	\item[$11.2$]	$\alpha_1=t$, $2\leq t\leq n-1$. If $\alpha_2\leq t-1$ then 
		$\alpha$ contains $213$, and if $t+1\leq\alpha_2\leq n-1$ then 
		$\alpha$ contains $(t,\alpha_2,1,n)$ or $(t,\alpha_2,n,1)$ which is 
		order-isomorphic to $213$ or $2341$ respectively, so $\alpha_2=n$. 
		Since $\alpha$ avoids $312$ we have that 
		$\alpha=(t,n,n-1,\dots,t+1,t-1,\dots,1)$, and there are $n-2$ 
		permutations of this type.

	\item[$11.2$] 	$\alpha_1=n$. Since $\alpha$ avoids $312$ we have that 
		$\alpha=(n,\dots,1)$.
\end{description}
Since the above cases are disjoint we obtain by corollary \ref{setword cor} 
that $|G_n|=|G_{n-1}|+n-2+1$. Besides $|G_4|=7$, hence similarly to the second 
proof we get $|G_n|={n\choose 2}+1$.\\

$12.$ By Simion and Schmidt ~\cite{sim} we have that $|S_n(213,321)|={n\choose 2}+1$, 
		hence by theorem $2.9$ we get $|S_n(213,321,\tau)|={n\choose 2}+1$ when 
		$\tau\in\{1324,2314\}$.\\ 

\textsc{Proof(13):} Let $A=\{231,312,1324\}$ and $G_n=S_n(r(A))=S_n(132,213,4231)$. By 
proposition \ref{symop} we get $|G_n|=|S_n(A)|$. Let $\alpha\in G_n$ and let us consider 
the possible value of $\alpha_1$:
\begin{itemize}
	\item[$13.1$]	$\alpha_1=1$. Since $\alpha$ avoids $132$ we have that 
		$\alpha=(1,2,\dots,n)$.

	\item[$13.2$]	$\alpha_1=t$, $2\leq t\leq n-1$. Similarly to case $6.1$ 
		we have that $\alpha=(t,t+1,\dots,n,\alpha_{n-t+2},\dots,\alpha_n)$.
		Evidently $\alpha\in G_n$ if and only if 
		$(\alpha_{n-t+2},\dots,\alpha_n)\in S_{t-1}(132,213,231)$. By Simion and 
		Schmidt ~\cite{sim} we have $|S_{t-1}(132,213,231)|=t-1$.

	\item[$13.3$]	$\alpha_1=n$. Evidently $\alpha\in G_n$ if and only if 
		$(\alpha_2,\dots,\alpha_n)\in S_{n-1}(132,213,231)$.
		By Simion and Schmidt ~\cite{sim} we have that $|S_{n-1}(132,213,231)|=n-1$.
\end{itemize}
Since the above cases are disjoint we obtain $|G_n|=1+(1+2+\dots+n-2)+n-1={n \choose 2}+1$.\\

$14.$ Let $A=\{213,132,4321\}$ and $G_n=S_n(r(A))=S_n(132,213,4231)$. By 
proposition \ref{symop} we get $|G_n|=|S_n(A)|$. Let $\alpha\in G_n$ and let us consider 
the possible value of $\alpha_1$:
\begin{itemize}
	\item[$14.1$]	$\alpha_1=1$. Since $\alpha$ avoids $132$ we have that 
		$\alpha=(1,2,\dots,n)$.

	\item[$14.2$]	$\alpha_1=t$, $2\leq t\leq n-2$. Fix $k$ such that $\alpha_k=n$, since
		$\alpha$ avoids $213$ and $132$ we have that $\alpha_i\geq t$ for $i\leq k$ 
		and $\alpha_i\leq t-1$ for $i\geq k+1$ respectively, so 
		$\alpha=(t,t+1,\dots,n,\alpha_{n-t+2},\dots,\alpha_n)$. Evidently 
		$\alpha\in G_n$ if and only if 
		$(\alpha_{n-t+2},\dots,\alpha_n)\in S_{t-1}(213,132,321)$. By Simion and 
		Schmidt ~\cite{sim} we have $|S_{t-1}(213,132,321)|=t-1$.

	\item[$14.3$]	$\alpha_1=n-1$. Since $\alpha$ avoids $213$ we have that $\alpha_2=n$. 
		Evidently $\alpha\in G_n$ if and only if 
		$(\alpha_3,\dots,\alpha_n)\in S_{n-2}(213,132,321)$.
		By Simion and Schmidt ~\cite{sim} we have that $|S_{n-2}(213,132,321)|=n-2$.

	\item[$14.4$]	$\alpha_1=n$. Evidently $\alpha\in G_n$ if and only if 
		$(\alpha_2,\dots,\alpha_n)\in S_{n-1}(213,132,321)$.
		By Simion and Schmidt ~\cite{sim} we have that $|S_{n-1}(213,132,321)|=n-1$.
\end{itemize}
Since the above cases are disjoint we obtain 
$|G_n|=1+(1+2+\dots+n-3)+n-2+n-1={n \choose 2}+1$.
\end{proof}
\\

All these results we summarize in the table $2$ (next page).
{ \footnotesize
\begin{table}[h]
    \begin{tabular}{|l|l|l|l|} \hline
						&			&						&			\\
	\emph{Representative $T\in\mathcal{T}$} &\emph{$|\mathcal{T}|$} & \emph{$|S_n(T)|$ for $T\in\mathcal{T}$}  	& \emph{Reference}	\\ \hline\hline
	
								&	&						&			\\
	$\{\alpha_1,\alpha_2,\tau\}$ when $\tau\in S_4$, $\tau$ contains &	&					&			\\
	$\alpha$ and $\overline{\{\alpha_1,\alpha_2\}}=\overline{\{123,132\}}$, & $160$ & $2^{n-1}$			& theorem \ref{the0n0} \\ \hline 

								&	&						&			\\
	$\overline{\{123,132,3412\}}$, $\overline{\{123,132,4231\}}$&	&						&			\\
	$\overline{\{123,213,3412\}}$, $\overline{\{123,213,4231\}}$& 	&  						&			\\ 
	$\overline{\{132,213,3412\}}$, $\overline{\{132,231,1234\}}$&	&						&			\\
	$\overline{\{132,231,2134\}}$, $\overline{\{132,231,3124\}}$& 	&  						&			\\ 
	$\overline{\{132,231,3214\}}$, $\overline{\{213,312,2341\}}$& $118$ & ${n\choose 2}+1$				& theorem \ref{thenn2}	\\
	$\overline{\{213,312,1324\}}$, $\overline{\{213,321,2314\}}$& 	&  						&			\\ 
	$\overline{\{231,312,1324\}}$, $\overline{\{132,213,4321\}}$&   &						&			\\
	$\{123,231,\tau\}$ when $\tau\in S_4$ contains 		&	&						&			\\
	$123$ or $231$						& 	&	  					&			\\ \hline

								&	&						&			\\
	$\{123,321,\tau\}$, $\tau\in S_4$ 			&	&						& Erd\"os and		\\
	$\{123,\alpha,4321\}$, $\alpha\in S_3$, $\alpha\neq 123$ & $32$	& $0$						& Szekeres ~\cite{erd} 	\\
	$\{321,\alpha,1234\}$, $\alpha\in S_3$, $\alpha\neq 321$ &	&						&			\\ \hline

								&	&						&			\\
	$\overline{\{123,312,\tau\}}$, when $\tau=1432$, 	&	&						&			\\
	$\tau=2143$, $2431$, $3214$, $3241$ or $3421$ 		& $24$	& $2n-2$					& theorem \ref{the2n2}  \\ \hline

								&	&						&			\\
	$\overline{\{123,132,3241\}}$, $\overline{\{132,213,2341\}}$& $12$ & $f_{n+2}-1$				& theorem \ref{thefn2}	\\ \hline

								&	&						&			\\
	$\overline{\{123,132,3421\}}$, $\overline{\{123,213,3421\}}$& $8$  & $3n-5$					& theorem \ref{the3n5}	\\ \hline

								&	&						&			\\
	$\overline{\{123,132,3214\}}$, $\overline{\{123,213,1432\}}$&   & $a_n$, where $a_n$ is the $n$-th  		&	 		\\ 
	$\overline{\{132,213,1234\}}$				& $6$  	& Tribonacci number ~\cite{fei} 		& theorem \ref{thtri}	\\ \hline	
    \end{tabular}
\caption{\small Cardinalities of the sets $S_n(\alpha_1,\alpha_2,\tau)$ 
		   when $\alpha_1,\alpha_2\in S_3$ and $\tau\in S_4$.}
\label{tab334}
\end{table} }
\section{$|S_n(T,\tau)|$ when $T\subset S_3$, $|T|=3$ and $\tau\in S_4$}
In this section we calculate the cardinalities of all the sets $S_n(T,\tau)$ when
$|T|=3$, $T\subset S_3$ and $\tau$ is any permutation in $S_4$.
\begin{theo}
\label{thefn}
  Let $T=\{123,132,213\}$ and let $\tau\in S_4$ contain 
  at least one permutation in $T$. Then $|S_n(T,\tau)|=|S_n(T)|=f_{n+1}$, 
  where $f_n$ is the $n$-th Fibonacci number.
\end{theo}
\begin{proof} By Simion and Schmidt ~\cite{sim} we have that $|S_n(T)|=f_{n+1}$, hence
by theorem \ref{opmu theo} we get $|S_n(T,\tau)|=f_{n+1}$.
\end{proof}

\begin{theo}
\label{the0}
   Let $T\subset S_3$ and $|T|=3$. For all $n\geq 6$,
				$$|S_n(T,\tau)|=0,$$
in the following cases:
\begin{enumerate}
	\item	$123\in T$ and $\tau=4321$.

	\item	$321\in T$ and $\tau=1234$. 

	\item	$\{123,321\}\subset T$ and $\tau\in S_4$.
\end{enumerate} 
\end{theo}
\begin{proof} By Erd\"os and Szekeres ~\cite{erd}. 
\end{proof}

\begin{theo}
\label{then}
  For all $n\in\mathcal{N}$,
			$$|S_n(T,\tau)|=n,$$
in the following cases:
\begin{enumerate}
	\item	$T=\{123,132,231\}$ and $\tau\in S_4$ contains at least one permutation 
		in $T$.

	\item	$T=\{123,132,213\}$ and $\tau=3412$.
\end{enumerate} 
\end{theo}
\begin{proof} $1.$ By Simion and Schmidt ~\cite{sim} we have that $|S_n(T)|=n$, hence by theorem 
\ref{opmu theo} we get $|S_n(T,\tau)|=n$ when $\tau\in S_4$ contains a permutation in $T$.\\

$2.$ Let $\alpha\in G_n=S_n(123,132,213,3412)$ and let us consider the 
possible value of $\alpha_1$:
\begin{itemize}
	\item[$2.1$]	$\alpha_1\leq n-2$. Since $\alpha$ avoids $123$ we have that 
		$\alpha$ contains $(\alpha_1,n,n-1)$, which means that $\alpha$ contains 
		$132$, a contradiction.

	\item[$2.2$]	$\alpha_1=n-1$. Since $\alpha$ avoids $213$ we get $\alpha_2=n$. 
		Since $\alpha$ avoids $3412$ we have only one permutation 
		$(n-1,n,n-2,\dots,1)$.

	\item[$2.3$]	$\alpha_1=n$. Evidently $\alpha\in G_n$ if and only if 
	$(\alpha_2,\dots,\alpha_n)\in G_{n-1}$.
\end{itemize}
Since the above cases are disjoint we obtain $|G_n|=|G_{n-1}|+1$. Besides $|G_4|=4$,
hence $|G_n|=n$ for all $n\geq 4$. It is easy to see that for $n=1,2,3$ the same
formula holds.
\end{proof}

\begin{prop}
\label{fromnton1}
  Let $T\subseteq S_k$. If $\alpha\notin S_n(T)$ then
  $(\alpha_1,\dots,\alpha_{j-1},n+1,\alpha_j,\dots,\alpha_n)\notin S_{n+1}(T)$ for 
  all $1\leq j\leq n$.
\end{prop}
\begin{proof} By definitions. 
\end{proof}

\begin{theo}
\label{the4}
  For all $n\geq 4$,{\small
 \begin{enumerate}
	\item	$S_n(123,132,213,3421)=\{(n-1,n,n-2,\dots,1)$, $(n-1,n,n-1,\dots,1,2)$, 
		$(n,\dots,3,1,2)$, $(n,\dots,1)\}$.

	\item	$S_n(123,132,213,4231)=\{(n,\dots,5,3,4,1,2)$, $(n,\dots,4,2,3,1)$, 
		$(n,\dots,3,1,2)$, $(n,\dots,1)\}$.
\end{enumerate} }
\end{theo}
\begin{proof} By induction and proposition \ref{fromnton1} .
\end{proof}

\begin{theo}
\label{the3}
  $\delta_n=(1,2,\dots,n)$. For all $3\leq n$, {\small
\begin{enumerate}
	\item	$S_n(123,132,231,3214)=\{(n,\dots,4,2,1,3)$, $(n,\dots,3,1,2)$, 
		$\delta_n\}$.

	\item	$S_n(123,132,231,4312)=\{(n-1,\dots,1,n)$, $(n,n-2,\dots,1,n-1)$, 
		$\delta_n\}$.

	\item	$S_n(123,132,231,4213)=\{(n-1,\dots,1,n)$, $(n,\dots,3,1,2)$, 
		$\delta_n\}$.

	\item	$S_n(123,231,312,1432)=\{(n-2,\dots,1,n,n-1)$, $(n-1,\dots,1,n)$, 
		$\delta_n\}$.

	\item	$S_n(123,231,312,2143)=\{(2,1,n,\dots,3)$, $(n-1,\dots,1,n)$, 
		$\delta_n\}$.

	\item	$S_n(132,213,231,1234)=\{(n,\dots,4,1,2,3)$, $(n,\dots,3,1,2)$, $\delta_n\}$.

	\item	$S_n(132,213,231,4123)=\{r(\delta_n)$, $(n,\dots,3,1,2)$, $\delta_n\}$.

	\item	$S_n(132,213,231,4312)=\{r(\delta_n)$, $(n,1,\dots,n-1)$, 
		$\delta_n\}$.

	\item	$S_n(132,213,231,4321)=\{r(\delta_n)$, $(n,1,\dots,n-1)$, 
		$(n,n-1,1,\dots,n-2)\}$.
\end{enumerate} }  
\end{theo}
\begin{proof} By induction and Proposition \ref{fromnton1}.
\end{proof}
\\

All these results we summarize in the following table.
{ \footnotesize
\begin{table}[h]
    \begin{tabular}{|l|l|l|l|} \hline
						&			&					&			\\
	\emph{Representative $T\in\mathcal{T}$}	&\emph{$|\mathcal{T}|$} &  \emph{$|S_n(T)|$ for $T\in\mathcal{T}$} & \emph{Reference}	\\ \hline\hline

								&	&			&			\\
	$T\cup\{\tau\}$ when $\overline T=\overline{\{123,132,231\}}$&	&			&			\\
	and, $\tau$ contains one permutation in $T$ or $\tau=3412$& $282$ & $n$			& theorem \ref{then} 	\\ \hline 

								&	&			&			\\
	$T\cup\{\tau\}$ when $123,321\in T$ or 			& 	&			& Erd\"os and		\\
	($123\in T$ and $\tau=4321$) or ($321\in T$ and $\tau=1234$)& $108$ & $0$		& Szekeres ~\cite{erd}	\\ \hline

								&	&			&			\\
	$\overline{\{123,132,231,\tau\}}$, $\tau=3214$, $4312$ or $4213$ &	&		&			\\
	$\overline{\{123,213,231,\tau\}}$, $\tau=1432$, $4132$ or $4312$ &	&		&			\\
	$\overline{\{123,231,312,\tau\}}$, $\tau=1432$, $2143$ or $3214$ & $46$	& $3$		& theorem \ref{the3}	\\
	$\overline{\{132,213,231,\tau\}}$, $\tau=1234$, $4123$, $4321$ or $4312$ &	&	&			\\ \hline

								&	&			&			\\
	$T\cup\{\tau\}$ when $\overline T=\{123,132,213\}$	&	&			&			\\
	and, $\tau$ contains one permutation in $T$  		& $38$ 	& $f_{n+1}$		& theorem \ref{thefn} 	\\ \hline 

								&	&			&			\\
	$\overline{\{123,231,312,\tau\}}$, $\tau=3421$ or $4231$& $6$	& $4$			& theorem \ref{the4}	\\ \hline
    \end{tabular}
\caption {\small Cardinalities of the sets $S_n(T,\tau)$ when 
			$T\subset S_3$, $|T|=3$ and $\tau\in S_4$.}
\label{tab3334}
\end{table} }

\section{$|S_n(T,\tau)|$ when $T\subset S_3$, $|T|=4,5,6$ and $\tau\in S_4$}
In this section we calculate the cardinalities of all the sets $S_n(T,\tau)$ when
$|T|\geq 4$, $T\subset S_3$ and $\tau$ is any permutation in $S_4$.
\begin{theo}
\label{the1}
  For all $n\geq 4$,{\small
	$S_n(123,132,213,231,4312)=$ $S_n(123,132,231,312,3214)=$ 
	$S_n(132,213,231,312,1234)=\{(n,\dots,1)\}$. } 
\end{theo}
\begin{proof} By induction and proposition \ref{fromnton1} .
\end{proof}

\begin{theo}
\label{the2b}
 Let $T\subset S_3$, $|T|=4$, $\{123,321\}\not\subset T$ 
  and let $\tau\in S_4$ contains at least one permutation in $T$. Then 
$|S_n(T,\tau)|=2$.
\end{theo}
\begin{proof} By Simion and Schmidt ~\cite{sim} we have that $|S_n(T)|=2$, so by theorem 
\ref{opmu theo}.
\end{proof}

\begin{theo}
\label{the1b}
 Let $T\subset S_3$ and $|T|=5$. Then for all $n\geq 3$, 
			$$|S_n(T,\tau)|=1,$$ 
in the following cases:
\begin{enumerate}
	\item	$S_n(T,\tau)=\{(1,2,\dots,n)\}$ if $321\notin T$ and $\tau\neq 4321$.

	\item	$S_n(T,\tau)=\{(n,n-1,\dots,1)\}$ if $123\notin T$ and $\tau\neq 1234$.
\end{enumerate}
\end{theo}
\begin{proof} By induction and proposition \ref{fromnton1} .
\end{proof}

{ \footnotesize
\begin{table}[h]
    \begin{tabular}{|l|l|l|l|} \hline
						&			&			        	&		   	\\
	\emph{Representative $T\in\mathcal{T}$}	&\emph{$|\mathcal{T}|$} & \emph{$|S_n(T)|$ for $T\in\mathcal{T}$}& \emph{Reference} 	\\ \hline\hline  

								&	&			&			\\
	$T\cup\{\tau\}$, $|T|\geq 4$ when $123,321\in T$, or 	& 	&			& Erd\"os and		\\
	$123\in T$ and $\tau=4321$, or $321\in T$ and $\tau=1234$& $348$ & $0$			& Szekeres ~\cite{erd}	\\ \hline

								&	&			&			\\
	$T\cup\{\tau\}$, $|T|=4$, $\{123,321\}\not\subset T$	&	&			&			\\
	and $\tau\in S_4$ contains permutation in $T$		& $100$ & $2$			& theorem \ref{the2b}	\\ \hline

								&	&			&			\\
	$T\cup\{\tau\}$ when					&	&			&			\\
	$|T|=5$, $123\notin T$ and $\tau\neq 1234$, or		&	&			&			\\
	$|T|=5$, $321\notin T$ and $\tau\neq 4321$, or		&	&			&			\\
	$\overline{\{123,132,213,231,4312\}}$,                	& $56$	& $1$			& theorem \ref{the1} \ref{the1b} \\
	$\overline{\{123,132,231,312,3214\}}$,                	&	&			&			\\
	$\overline{\{123,213,231,312,1432\}}$,                	&	&			&			\\
	$\overline{\{132,213,231,312,1234\}}$,                	&	&			&			\\ \hline

    \end{tabular}
\caption {\small Cardinalities of the sets $S_n(T,\tau)$ when $T\subset S_3$, $|T|\geq 4$ and $\tau\in S_4$.}
\label{tab33334}
\end{table} }

\begin{theo}
\label{the0b}
   Let $T\subset S_3$. Then for all $n\geq 6$,
				$$|S_n(T,\tau)|=0,$$
in the following cases:
\begin{enumerate}
	\item	$123\in T$ and $\tau=4321$.

	\item	$321\in T$ and $\tau=1234$. 

	\item	$\{123,321\}\subset T$ and $\tau\in S_4$.
\end{enumerate} 
\end{theo}
\begin{proof} By Erd\"os and Szekeres ~\cite{erd}. 
\end{proof}
\\

All these results we summarize in the table [3].

\end{document}